\documentclass[11pt]{amsart}
\usepackage[utf8]{inputenc}

\usepackage[nospace,noadjust]{cite}
\usepackage{relsize}
\usepackage{cite}
\usepackage{color}
\usepackage[dvipsnames]{xcolor}

\usepackage{tikz-cd}
\usetikzlibrary{calc,fit,matrix,arrows,automata,positioning}
\usetikzlibrary{decorations.pathreplacing,angles,quotes}
\usepackage{colortbl}
\usepackage{hyperref, enumitem}
\hypersetup{
  colorlinks   = true, 
  urlcolor     = blue, 
  linkcolor    = Purple, 
  citecolor   = red 
}
\usepackage{amsmath,amsthm,amssymb}

\usepackage[margin=2.5cm]{geometry}

\usepackage{stmaryrd}

\usepackage{array}
\newcolumntype{x}[1]{>{\centering\arraybackslash\hspace{0pt}}p{#1}}

\theoremstyle{definition}
\newtheorem{theorem}{Theorem}[section]
\newtheorem{definition}[theorem]{{{Definition}}}

\newtheorem{proposition}[theorem]{{{Proposition}}}
\newtheorem{lemma}[theorem]{{{Lemma}}}

\newcommand{\numberset}{\mathbb}

\newcommand{\C}{\mathcal{C}}

\newcommand{\F}{\numberset{F}}

\newcommand{\mG}{\mathcal{G}}
\newcommand{\mL}{\mathcal{L}}

\newcommand{\mI}{\mathcal{I}}

\newcommand{\mM}{\mathcal{M}}

\newcommand{\mU}{\mathcal{U}}

\renewcommand{\mG}{\mathcal{G}}
\newcommand{\mV}{\mathcal{V}}

\newcommand{\mH}{\mathcal{H}}

\newcommand{\wt}{\textnormal{wt}}

\newcommand{\Fq}{\F_q}
\newcommand{\Fm}{\F_{q^m}}

\newcommand{\Fqn}{\mathbb{F}_{q}^{n}}
\newcommand{\Fqm}{\mathbb{F}_{q^{m}}}
\newcommand{\Fqmn}{\mathbb{F}_{q^{m}}^{n}}

\newcommand{\rk}{\textnormal{rk}}

\newcommand{\D}{\textnormal{D}}

\DeclareMathOperator{\PG}{PG}

\newcommand{\rowspan}{\textnormal{rowspan}}

\title{Recent advances on minimal codes}
\usepackage[foot]{amsaddr}
\author{Martin Scotti$^1$}
\address{$^1$Universit\'e Paris 8, Laboratoire de G\'eom\'etrie, Analyse et Applications, LAGA, Universit\'e Sorbonne Paris Nord, CNRS, UMR 7539, France.}

\email{martin.scotti@etud.univ-paris8.fr}
\thanks{M.~S. is partially supported by the ANR-21-CE39-0009 - BARRACUDA (French \emph{Agence Nationale de la Recherche}).}

\begin{document}

\begin{abstract} 
In this short survey we concern ourselves with minimal codes, a classical object in coding theory. We will explain the relation between minimal codes and various other mathematical domains, in particular with finite projective geometry. This latter connection has sparked a renewed interest in minimal codes, giving rise to new constructions as well as new questions.
\end{abstract}

\maketitle


\bigskip

\section{Introduction}

Minimal codes are a classical object in coding theory, and their simple definition provides natural connections with various areas of mathematics. Minimal codewords were first introduced by Hwang in \cite{hwang}, and attracted renewed interest after they were studied by Massey in \cite{Massey} in order to construct an efficient secret sharing scheme. To the best of our knowledge, minimal codes, i.e. codes where all nonzero codewords are minimal, are first studied by Ashikhmin and Barg in \cite{ashikhmin1998minimal}. Research on minimal \emph{codewords} in various contexts is an active area of research (see \cite{Bartoli2022MinimalCI} for a recent example), but in the present survey we shall only examine minimal codes. Historically, research has first focused on the theory of binary minimal codes, where the definition of a minimal code coincides with the definition of an intersecting code. A rich theory has developed, with several natural generalizations and applications, see for example \cite{cohen1994intersecting, CL, mesnager2019several}.

\medskip

The main reason for renewed interest in minimal codes in recent years has been a new geometric interpretation of minimal codes, first outlined in 2019 simulateously in \cite{ABNgeo} and \cite{tang}. As explained in these works, the geometric counterpart of a minimal code is a set of points in the projective space called a strong blocking set. These sets had been introduced by Davydov, Giulietti, Marcugini and Pambianco in \cite{davydov2011linear} for their connection with the theory of saturating sets. This geometric interpretation has led to a deeper understanding of the properties of minimal codes, as well as improved bounds on their parameters. Additionally, many explicit geometric constructions of strong blocking sets of small size have been obtained, yielding explicit constructions of minimal codes with small length. It should be noted that efficient explicit constructions of minimal codes with small parameters are not only of theoretical interest, but also allow for more efficient applications of minimal codes in practice, as we shall see at the end of the paper.

\medskip

The objective of this survey is to present these recent advances as well as some classical applications of the theory of minimal codes to various other areas of mathematics. We hope that this work will clarify which questions remain open for future research and which techniques are used to obtain new geometric insights.

\medskip

The survey is organized as follows. In Section \ref{sec:def} we fix some notation and provide an elementary definition of minimality that is independent of the metric used, as well as a definition of a strong blocking set. In Section \ref{sec:hamming} we investigate the case of minimal codes in the Hamming metric. Since this imposes no additional restrictions on the corresponding strong blocking sets, many results will take the form of general bounds and constructions in the theory of strong blocking sets. In Section ~\ref{sec:rank} we examine minimal codes in the rank metric. In that case, is the corresponding linear sets that must be strong blocking sets, which gives rise to a different set of bounds and constructions. Finally, in Section \ref{sec:open} we present some open problems in the theory of minimal codes for future research as well as some applications and connections with other mathematical problems.

\section{Background}\label{sec:def}

\subsection{Coding theory}

We begin by introducing classical notation from coding theory and projective geometry.

Let $q$ be a prime power, and $\Fq$ the finite fiels with $q$ elements. In the Hamming metric, a code $\C$ will be a vector subspace of $\Fqn$, while in the rank metric, it will be a subspace of $\Fqmn$. A matrix $G$ such that $\C$ is equal to the rowspan of $G$ is called a generator matrix. The code is called \emph{degenerate} if a generator matrix has a zero column, and \emph{nondegenerate} if it is not degenerate. We will endow $\Fqn$ with a notion of support, defined as follows.

\begin{definition}
Let $P$ be a poset. A support is an application $\sigma : \Fqn \rightarrow P$ such that
$$\forall x \in \Fqn\quad \forall \lambda \in \Fq^{\times} \quad \sigma(\lambda \cdot x) = \sigma(x),$$
and
$$\forall x \in \Fqn \quad \sigma(0) \leq \sigma(x).$$
\end{definition}

Now we can define what a minimal code is in a very general way.

\begin{definition}
Let $\C \subset \Fqn$ be a code. A codeword $c \in \C$ is said to be \emph{minimal} if there is no $c' \in \C\setminus\{0\}$ such that $\sigma(c') \subsetneq \sigma(c)$.
The whole code $\C$ is said to be \emph{minimal} if each nonzero codeword $c \in \C\setminus\{0\}$ is minimal.
\end{definition}

An equivalent definition is that a code is minimal if and only if the set of its nonzero supports forms an antichain.

\subsection{Finite geometry}

We define the projective space $\PG(k-1, q)$ as the set of equivalence classes of points of $\Fq^{k} \setminus \{0\}$ with respect to the collinearity equivalence relation. It is possible to map any point of $\Fq^{k} \setminus \{0\}$ to its equivalence class, which is a point of $\PG(k-1, q)$. A hyperplane of $\PG(k-1, q)$ is the image of a subspace of codimension $1$ of $\Fqn$ under this map.

\begin{definition}
A strong blocking set (first defined in \cite{davydov2011linear}) is a set of projective points $S \subset \PG(k-1, q)$ such that its intersection with any projective hyperplane $\mH$ spans that hyperplane. Equivalently, it is a set of points capable of generating any projective hyperplane.
\end{definition}

A \emph{linear set} of $\PG(k-1, q^{m})$ (introduced in \cite{Lunardon1999NormalS}) is the image of a $\Fq$-linear vector subspace of $\Fm^{k}$ under the projection map.

\section{The Hamming metric case} \label{sec:hamming}

In this section we first provide some definitions of codes in the Hamming metric and of minimal codes in this context. After this we shall move on to their properties as well as their geometric counterpart: strong blocking sets.

\subsection{Definitions}

The support of a vector $x \in \Fqn$ is defined as
$$\sigma(x) = \{1\leq i\leq n\mid x_{i} \neq 0\} \subseteq \{1, \dots, n\}$$
and the Hamming weight of $x$ is then defined to be $\wt(x) = |\sigma(x)|$. The Hamming weight induces a distance called the Hamming metric: for two vectors $x, y \in \Fq$ one has $d(x, y) = \wt(y-x)$.

We consider linear codes, that is vector subspaces of $\Fqn$.
For such a code $\C$ it is customary to write its dimension $k = \dim(\C)$, and its minimum distance $d = d(\C) = \min_{c \in \C \setminus \{0\}} \wt(c)$.
The code is then said to have parameters $[n, k, d]_{q}$.

\medskip

With this definition of support, the definition of minimal code in the context of the Hamming metric goes as follows. A code $\C$ is a minimal code if for every two nonzero codewords $c, c' \in \C\setminus \{0\}$ one has
$$\sigma(c) \subseteq \sigma(c') \iff \exists \lambda \in \Fq, \quad c = \lambda c'.$$

\smallskip

A code $\C$ is said to be intersecting if for any two nonzero codewords $c, c' \in \C\setminus\{0\}$, one has
$$\sigma(c) \cap \sigma(c') \neq \varnothing.$$
It is straightforward to check that when $q=2$, $\C$ is an intersecting code if and only if it is minimal. In general, minimal codes are intersecting codes, but the reverse is not true if $q > 2$.

\subsection{Bounds on the parameters of Hamming-metric codes}

We begin with an upper bound on the rate of minimal codes.

\begin{proposition}[\cite{chabanne2013towards}, Theorem 2.]
Let $\C$ be a minimal $[n, k, d]_{q}$-code. Then
$$k/n \leq \log_{q}(2).$$
\end{proposition}

This proposition implies in particular that minimal codes must be $\emph{rather long}$ (i.e. that $n$ cannot be too small compared with $k$. The geometric interpretation of minimal codes exposed below will provide a much more detailed picture.

The following proposition establishes a bound on the weight of every codeword analogous to the Singleton bound.

\begin{proposition}[\cite{ashikhmin1998minimal}, Lemma 2.1. (2.)]
Let $\C$ be a minimal $[n, k, d]_{q}$-code. Then for any codeword $c\in \C$ one has
$$\wt(c) \leq n-k+1.$$
\end{proposition}

A direct consequence of these two propositions is that the weight distribution of a minimal code must be rather restricted. Most importantly, the maximum weight of a codeword cannot be too high. It is possible to establish a sufficient condition concerning the weight distribution of a code for it to be minimal, as done in \cite{ashikhmin1998minimal}. It is now known as the Ashikhmin-Barg condition.

\begin{theorem}[\cite{ashikhmin1998minimal}, Lemma 2.1. (3.)]
Let $\C$ be a $[n, k, d]_{q}$-code, and let $d \leq w \leq n$ be its maximal Hamming weight.
Assume that
$$\frac{w}{d} < \frac{q}{q-1}.$$
Then $\C$ must be a minimal code.
\end{theorem}

As an immediate consequence, one-weight codes must be minimal codes. This is used in \cite{mesnager2019several} in order to construct minimal codes with few weights. It must also be noted that the Ashikhmin-Barg condition is not necessary, in fact there exists an infinite family of minimal codes violating the condition, as established in \cite{HENG2018176}.

These preliminary considerations, while interesting, miss part of the picture. A substantial leap is accomplished by examining the geometric counterpart of minimal codes: strong blocking sets.

\subsection{Strong blocking sets} \label{subs:sbs}

The geometric idea, explored simultaneously and independently in \cite{ABNgeo} and \cite{tang2019full} is to take the columns of a generator matrix $G$ of a nondegenerate minimal code $\C$ with parameters $[n, k, d]_{q}$ and to consider the corresponding points in the projective space $\PG(k-1, q)$.
This is a classical procedure in coding theory, and it yields a multiset (i.e. a set with some repetitions, as $G$ may contain collinear columnns) of size $n$.

The crucial observation is the following.

\begin{theorem}[\cite{ABNgeo}, Theorem 3.4., \cite{tang2019full}, Theorem 14.]
$\C$ is a minimal code if and only if the pointset obtained from any generator matrix is a strong blocking set.
\end{theorem}

This means that studying strong blocking sets is equivalent to studying minimal codes.

Therefore, it is possible to use many geometric techniques to study minimal codes.
In particular, it is clear from their definition that if $S \subset \PG(k-1, q)$ is a strong blocking set and $S \subseteq S'$, then $S'$ must also be a strong blocking set. Therefore, a very natural question is that of determining the smallest size of a strong blocking set in $\PG(k-1, q)$. This smallest size is denoted by the function $m(k, q)$, it corresponds both to the shortest length of a minimal code of dimension $k$ over $\Fq$, and to the smallest size of a strong blocking set in $\PG(k-1, q)$. While it may not have been clear from the definition of minimal codes that they must be \emph{rather long}, it is very clear from the geometric point of view that short minimal codes (and any bounds on the function $m(k,q)$) are of great interest.

\begin{proposition}[\cite{ABNgeo}, Theorem 4.3.]
Let $\C$ be a minimal $[n, k, d]_{q}$-code with $k \geq 2$. Then
$$k \leq d + q - 2.$$
\end{proposition}

\begin{proposition} [\cite{3CB}, Theorem 2.8.]
Let $\C$ be a minimal $[n, k, d]_{q}$-code. Then
$$d \geq (q-1)(k-1) + 1.$$
\end{proposition}

This, as well as following lower bound, are established by Alfarano, Borello, Neri and Ravagnani in \cite{3CB}.

\begin{theorem}[\cite{3CB}, Theorem 2.14.]\label{lower bound}
$$m(k, q) \geq (q+1)(k-1).$$
\end{theorem}

This is improved slightly in \cite{scotti2024lower} and \cite{Bishnoi_trifference}, where the following asymptotic improvement is given.

\begin{theorem}[\cite{scotti2024lower}, Theorem 3.3., and \cite{Bishnoi_trifference}, Theorem 1.4.]
$$\liminf_{k \rightarrow \infty} \frac{m(k, q)}{k} \geq q+ \varepsilon(q),$$
where $\varepsilon$ is an increasing function such that $1.52 <\varepsilon(2)$ and $\lim_{q \rightarrow \infty} \varepsilon(q) = \sqrt{2} + \frac{1}{2}$.
\end{theorem}

A natural consequence of this is that there are only finitely many minimal codes reaching the lower bound of Theorem \ref{lower bound} for any fixed $q$.

These bounds can be interpreted as impossibility results, i.e. strong blocking sets cannot be smaller than $m(k, q)$.

Conversely, upper bounds on $m(k, q)$ represent existence results on strong blocking sets. These can either be nonconstructive (for instance with probabilistic arguments), or stem from an explicit construction.
These are best understood in their geometric statement regarding strong blocking sets.

For instance, if one takes points at random in $\PG(k-1, q)$, at some point the probability that their union is a strong blocking set is positive. This means that there must be a strong blocking set of that size (while not providing an explicit example of such a strong blocking set). This approach is considered in \cite{cohen1994intersecting}, and yields the following bound.

\begin{theorem}[\cite{cohen1994intersecting}, Theorem 8.1.]
$$m(k, q) \leq \frac{2k}{\log_{q}(\frac{q^{2}}{q^{2}-q+1})} \simeq 2q\ln(q)k.$$
\end{theorem}

A further refinement of this elementary idea consists in taking not random points, but random lines. Indeed, since any line and any hyperplane of $\PG(k-1, q)$ must intersect, it is natural that sets that are unions of lines would have large intersection with any hyperplane, making them more suitable candidates for strong blocking sets.

This approach is considered in \cite{ABN2023} and \cite{Bishnoi_trifference}, where the authors use different methods to establish the same further improvement (a more careful consideration of the ideas presented in \cite{heger2021short} leads to the same improvement as well).

\begin{theorem}[\cite{ABN2023}, Theorem 4.10. and \cite{Bishnoi_trifference}, Theorem 1.1.]
$$m(k, q) \leq \left\lceil \frac{2k}{\log_{q}(\frac{q^{4}}{q^{3} - q + 1})}\right\rceil  \cdot (q+1) \simeq 2(q+1)k.$$
\end{theorem}

\subsection{Explicit constructions of strong blocking sets}

A first construction of a strong blocking set is the tetrahedron, described for example in \cite{ABNgeo}.
Consider $k$ points $P_{1}, \dots, P_{k}$ in general position. For each pair of points $P_{i}, P_{j}$, we call $E_{ij}$ the line passing through $P_{i}$ and $P_{j}$. Define $\mL$ as the union of all points on the lines $E_{ij}$. Then $\mL$ is a strong blocking set of size $(q+1) {{k}\choose{2}} - k$.

\medskip

A common tool for the construction of minimal codes is concatenation, defined in following way.

\begin{definition}
Let $\C$ be an $[N, K, D]_{q^{k}}$-code and $\mI$ an $[n, k, d]_{q}$-code. Define $\phi : \Fq^{k} \rightarrow \Fqn$ a linear map with image $\mI$.
Then the concatenation of $\C$ and $\mI$ by $\phi$ is
$$\mI \square_{\phi} \C = \{(\phi(c_{1}, \dots, \phi(c_{n}))\mid c\in \C\}.$$
\end{definition}

If one imposes restrictions on the codes $\mI$ and $\C$ it is possible to guarantee that the concatenation will be a minimal code. This allows the construction of asymptotically good minimal codes, where one chooses a sequence of AG codes over a large field $\mathbb{F}_{q^{k}}$, together with a well-chosen minimal inner code $\mI$. Such an approach was outlined by Cohen, Mesnager and Patey in \cite{cohen2013minimal} and fully developped by Bartoli and Borello in \cite{bartoli2021small}, yielding an explicit construction of strong blocking sets with size $O(kq^{4})$.

\medskip

A more direct investigation into AG codes done by Randriambololona in \cite{randriambololona20132} yields efficient constructions of AG codes that are directly constructed to be intersecting. Together with the concept of outer minimal codes developped by Alfarano, Borello and Neri in \cite{ABN2023}, they allow an efficient concatenation construction in the binary case, yielding an explicit construction of minimal codes with asymptotic rate $6/35$, as noted in \cite{BSSintersecting}.

\medskip

In \cite{Fancsali}, Fancsali and Sziklai introduce the avoidance property (or sets of lines in higgledy-piggledy arrangement), defined as follows.

\begin{definition}
Let $\mL$ be a set of lines of $\PG(k-1, q)$. This set is said to have \emph{avoidance property} if for any projective subspace $\mV$ of codimension $2$, there exists a line $\ell \in \mL$ such that $\ell \cap \mV = \varnothing.$
\end{definition}

It turns out that any set of lines with avoidance property is a strong blocking set. Therefore, in order to give an explicit construction of strong blocking sets, one possibility is to try to construct an explicit set of lines with avoidance property.
In \cite{fancsali2016higgledy}, the authors give bounds on the minimum size that a set with avoidance property must have.

\begin{theorem}[\cite{Fancsali}, Theorem 14. and Theorem 20.]
Let $\mL \subset \PG(k-1, q)$ be a set of lines with avoidance property. Assume that $q \geq \lfloor (k-1)/2\rfloor + (k-1)$. Then $\mL$ must contain at least $\lfloor (k-1)/2\rfloor + (k-1)$ lines.

Furthermore, if $q \geq 2k-3$, then there exists an explicit set of $2k-3$ lines (tangents to the rational normal curve) with avoidance property.
\end{theorem}

\medskip

The most efficient construction known so far has been obtained by Alon, Bishnoi, Das and Neri in \cite{SBSfromexpander}, and it uses a graph-theoretic property to construct sets of lines with avoidance property.

Let $G = (V, E)$ be a finite connected graph. Its vertex integrity is defined as
$$\iota(G) = \max_{S \subseteq V} |S| + \kappa(G-S),$$
where $\kappa(H)$ is the size of the largest connected component of $H$.

Consider points $P_{1}, \dots, P_{n} \in \PG(k-1, q)$ that form a projective $[n, k, d]_{q}$-system (i.e. the projective counterpart of an $[n, k, d]_{q}$-code as defined in Subsection \ref{subs:sbs}), and let $\mG = (V, E)$ be a graph on $n$ vertices. For each pair of points $P_{i}, P_{j}$, we call $\ell_{i, j}$ the line passing through both of them. Let
$$\mL = \bigcup_{(E_{i}, E_{j})\in V} \ell_{i, j}.$$
Now the crucial insight is the following.

\begin{lemma}[\cite{SBSfromexpander}, Lemma 4.4.] \label{lem:exp}
Let $\mM$ be a projective $[n, k, d]_{q}$-system, and let $G = (\mM, E)$ be a graph with vertex integrity $\iota(G) \geq n-d+1$. Then $\mL$ is a set of lines with avoidance property.
\end{lemma}

The authors then use explicit constructions of algebraic geometry codes with suitable parameters, as well as explicit constructions of expander graphs with suitable vertex integrity, together with Lemma \ref{lem:exp} in order to obtain an explicit construction of strong blocking sets with small size. This yields the following theorem.

\begin{theorem}[\cite{SBSfromexpander}, Theorem 4.5.]
There exists an explicit construction of strong blocking sets with size $O(kq)$.
\end{theorem}

The implicit constant in this result depends on $q$. In general, it is slightly lower when $q$ is a square than when it is an odd power of a prime. Furthermore, the constant grows smaller as $q$ tends to infinity, and approaches a limit value of around 22. The interested reader is invited to check the details (which are beyond the scope of the present survey) in \cite{SBSfromexpander}.

\subsection{Optimal constructions in small dimension}

Determining with precision the size of the smallest strong blocking is possible when $k$ is small, although this quickly becomes a difficult computational problem as $k$ grows large. For low values of $k$ it also becomes possible to study the shape of strong blocking sets of minimal size, as done by Smaldore in \cite{smaldore}.

\smallskip

It seems that the different small values of $m(k, q)$ have not been recorded in any one work, mostly because the choice of $q$ is usually done beforehand. Below, we attempt to give the most precise information available for small values of $k$ and $q$, be it the precise value, or the best known bounds. Many advances have been obtained recently by Kurz in \cite{kurz2023divisible} and \cite{kurztrifferent}, in particular when $q=2$ and $k \geq 7$. The specific bound $m(10, 2) \leq 30$ is established by Cohen and Zémor in \cite{cohen1994intersecting}. The bounds on $m(6, 3)$ given here are derived in \cite{Bishnoi_trifference}. The reader is invited to refer to these works for more details on the methods used to establish them.

\begin{table}[ht!]
\begin{tabular}{c|ccccccccccc}
$k$ & $2$ & $3$ & $4$ & $5$  & $6$  & $7$  & $8$  & $9$  & $10$                    & $11$                   & $12$                    \\ \hline
$m(k, 2)$ & $3$ & $6$ & $9$ & $13$ & $15$ & $20$ & $24$ & $26$ & $28 \leq \cdot \leq 30$ & $31\leq \cdot \leq 35$ & $33 \leq \cdot \leq 40$
\end{tabular}
\end{table}

\begin{table}[ht!]
\begin{tabular}{c|ccccc}
$k$       & $2$ & $3$ & $4$  & $5$  & $6$                     \\ \hline
$m(k, 3)$ & $3$ & $9$ & $14$ & $19$ & $22 \leq \cdot \leq 24$
\end{tabular}
\end{table}

\section{The rank metric case} \label{sec:rank}

In this section we examine the case of minimal codes in the rank metric. These have been introduced by Alfarano, Borello, Neri and Ravagnani in \cite{rankminimal}.

\subsection{Definitions}

A rank-metric code $\C$ is a $\Fqm$-linear subspace of $\Fqmn$. Its dimension is $k = \dim_{\Fqm}(\C)$. Let $\Gamma$ be an $\Fq$-basis of $\Fqm$. We can define an $m\times n$ matrix for every codeword $c\in \C$ in basis $\Gamma$. We define the support of a codeword $c \in \C$ as
$$\sigma_{\rk}(c) = \rowspan(c) \subseteq \Fqn.$$
It is a classical observation to note that $\sigma_{\rk}(c)$ does not depend on the choice of $\Gamma$.

We can define the rank weight of a codeword $c\in \C$ to be $\wt_{\rk}(c) = \dim_{\Fq}(\sigma_{\rk}(c))$. This induces a distance over $\Fqmn$, called the rank distance. Similarly to the Hamming metric case, a code $\C \subseteq \Fqmn$ of dimension $k$ and minimum distance $d = \min_{c\in \C\setminus\{0\}}\wt_{\rk}(c)$ is said to have parameters $[n , k, d]_{q^{m}/q}$.

The definition of support in the rank metric is in accordance with that given in Section \ref{sec:def}. Consequently, in the rank metric, a minimal code is a code for which for every two nonzero codewords $c, c' \in \C\setminus \{0\}$ one has
$$\sigma_{\rk}(c) \subseteq \sigma_{\rk}(c') \iff \exists \lambda \in \Fqm, \quad c = \lambda \cdot c'.$$

\subsection{The geometric interpretation}

Before discussing the geometric interpretation of minimal codes in the rank metric, we must first define how projective sets are defined from a generator matrix $G$. We begin by introducing $q$-systems.

\begin{definition}
An $[n, k]_{q^{m}/q}$ system is an $\Fq$-linear vector subspace $\mU\subseteq \Fqm^{k}$ such that $\dim_{\Fq}(\mU)=n$.
\end{definition}

Let $\mU$ be an $[n, k]_{q^{m}/q}$ system.
The linear set associated to $\mU$ is the set of points 
\[L_\mU:=\{\langle x \rangle_{\F_{q^m}} \mid  x\in \mU\setminus\{0\}\}\subseteq \PG(k-1, q^m),\]
where $\langle x \rangle_{\F_{q^m}}$ denotes the projective point corresponding to $x$.

Now if $G$ is the generator matrix of a code $\C$ over $\Fqm$, we can define its corresponding $[n, k]_{q^{m}/q}$ system $\mU$ to be the $\Fq$-linear columnspan of columns of $G$.

The geometric interpretation of minimality outlined in \cite{rankminimal} is the following.

\begin{theorem}
Let $\C$ be a rank-metric code with parameters $[n , k, d]_{q^{m}/q}$ and let $\mU$ be a corresponding $[n, k]_{q^{m}/q}$ system obtained from any generator matrix.
Then $\C$ is a minimal code if and only if $L_{\mU}$ is a strong blocking set.
\end{theorem}

Note that wether $L_{\mU}$ is a strong blocking set does not depend on the choice of a generator matrix. A linear $[n, k]_{q^{m}/q}$ system $\mU$ is called a \emph{linear cutting blocking set} if $L_{\mU}$ is a strong blocking set.

\smallskip

As one can see, it turns out that strong blocking sets are still the geometric objects corresponding to minimal codes in the rank metric!

\subsection{Bounds on the parameters of minimal codes in the rank metric}

Similarly to minimal codes in the Hamming metric, the geometric interpretation of minimal codes in the rank metric as the counterpart of strong blocking (linear) sets yields various properties, some of which we present here.

\begin{proposition}[\cite{rankminimal}, Corollary 5.9.]
Let $\C$ be a minimal $[n, k, d]_{q^{m}/q}$-code in the rank metric. Let $c \in \C$ be a codeword. Then
$$\wt_{\rk}(c) \leq n - k + 1.$$
\end{proposition}

This allows to prove the following theorem.

\begin{theorem}[\cite{rankminimal}, Corollary 5.10.]\label{ranklb}
Let $\C$ be a rank-metric $[n, k, d]_{q^{m}/q}$-code. If $\C$ is minimal, then
$$n \geq k + m - 1.$$
\end{theorem}

Along with this lower bound, the geometric interpretation also yields upper bounds on the length of minimal codes in the rank metric.

\begin{theorem}[\cite{rankminimal}, Proposition 6.2.]
Let $\C$ be a nondegenerate $[n, k]_{q^{m}/q}$-code with $n \geq (k-1)m + 1$. Then $\C$
is minimal.
\end{theorem}

\begin{theorem}[\cite{rankminimal}, Theorem 6.11.]\label{rankub}
Let $k,m \geq 2$. For any set of parameters $n, k, m, q$ such that $n \geq 2k + m - 2$, there exists an  $[n, k]_{q^{m}/q}$-code that is minimal in the rank metric.
\end{theorem}

For a thourough discussion of the above results and a more complete picture, we encourage the reader to refer to \cite{rankminimal}. In particuler, note that Theorem \ref{ranklb} and Theorem \ref{rankub} taken together determine the possible lengths of a rank-metric minimal code, leaving out $k-1$ values of $n$ for which it is not known wether minimal codes exist.

In the particular case when $k = 3$, several explicit constructions of codes of length $n = m+2$ have been obtained both in \cite{rankminimal} and in the subsequent \cite{LLMT}.

\section{Applications and open problems} \label{sec:open}

As explained in the introduction, minimal codes were first studied for their relation with secret-sharing. In this section we briefly present several other applications as well as connections with different areas of mathematics. We shall also explain what we consider to be the main questions still open for future research.

\subsection{Applications and related problems}

In \cite{obliviouscrepeau}, Brassard, Crépeau and Santha outline the use of binary minimal codes for oblivious transfer protocols.

Suppose two parties, called Alice and Bob, want to perform the following protocol: Alice holds two $k$ bit strings $x_{0}, x_{1}\in \mathbb{F}_{2}^{k}$, and Bob holds a bit $b \in \{0, 1\}$. Bob wants to know $x_{b}$ without revealing $b$ to Alice, and Alice wants Bob to gain information on at most one of $x_{0}$ and $x_{1}$. We assume that Alice and Bob have access to a protocol that performs oblivious transfer in the case when $k = 1$.

Given a $[n, k]_{2}$-code that is minimal, the authors provide a protocol for $k$-bit oblivious transfer that uses $1$-bit oblivious transfer exactly $n$ times. The smaller $n$ is, the faster the protocol will be. Therefore, in this context, it is clearly of interest to have access to an efficiently computable explicit construction with small length.

\bigskip

The notion of covering radius of a code is connected to the study of saturating sets, a classical topic in projective geometry. These are defined in the following manner.

\begin{definition}
Let $S \subset \PG(k-1, q)$. If $\rho$ is the smallest integer such that any point $Q \in \PG(k-1, q)\setminus S$ verifies $Q \in \langle P_{1}, \dots, P_{\rho+1} \rangle$ (where the $P_{i}$'s are points of $S$), then the set $S$ is said to be $\rho$-saturating.
\end{definition}

Strong blocking sets actually provide an example of saturating sets, as shown in \cite{davydov2011linear}.

\begin{theorem}[\cite{davydov2011linear}, Theorem 3.2.]
Any strong blocking set in a subgeometry $\PG(k-1, q)$ of $\PG(k-1, q^{k-1})$
is a $(k-2)$-saturating set in $\PG(k-1, q^{k-1})$.
\end{theorem}

\medskip

In the rank metric, the covering radius is also equivalent to rank saturating systems, introduced by Bonini, Borello and Byrne in \cite{BBBranksaturating}, and further developped in \cite{BBMSaturatingminimalrank}. As explained in \cite{BBBranksaturating}, linear cutting blocking sets provide examples of saturating systems.

\begin{theorem}[\cite{BBBranksaturating}, Theorem 4.6.]
Let $\mU$ be an $[n, k]_{q^{m}/q}$ system. If $\mU$ is a linear cutting blocking set, then it is a $(k-2)$-saturating set in $\PG(k-1, q^{m(k-1)})$.
\end{theorem}

Saturating systems in the sum-rank matric, together with their connection with minimal codes (also in the sum-rank metric), are explored in the recent work \cite{BBBsumranksaturating}.

\bigskip

Another connection concerns additive combinatorics. Consider $G$ a finite abelian group. Its Davenport constant, noted $\D(G)$, is the smallest integer $\ell$ such that for any $\ell$ elements $a_{1}, \dots, a_{\ell} \in G$ (possibly with repetitions), there must be some of them that sum to the group identity:
$$\sum_{k = 1}^{r} a_{i_{k}} = 0_{G},$$
with $1 \leq i_{1} < \dots < i_{r} \leq \ell$.
Similarly it is possible to define the $2$-way Davenport constant, noted $\D_{2}(G)$, by requiring that there be $2$ disjoint zero-sum subsequences.

The study of $\D_{2}(C_{2}^{r})$, where $C_{2}$ is the group with $2$ elements, is equivalent to the theory of minimal binary codes, as pointed out by Schmid and Plagne in \cite{SP}.

\bigskip

One more related concept is that of trifferent codes, examined for instance in \cite{Bishnoi_trifference}, \cite{kurztrifferent} and more recently in \cite{Bhandari2024ImprovedUB}. A trifferent code is a subset $\C \subseteq \mathbb{F}_{3}^{n}$ such that for any $3$ distinct elements $c_{1}, c_{2}, c_{3} \in \C$ there is a coordinate $i$ for which $\{c_{1, i}, c_{2, i}, c_{3, i}\} = \{0, 1, 2\}$.
A classical problem is to determine the largest possible trifferent code in $\mathbb{F}_{3}^{n}$. This quantity is noted $T(n)$ and has been studied intensively. To the best of our knowledge, the most recent improvement on explicit small values have been obtained by Kurz in \cite{kurztrifferent}, while the best asymptotic upper bound on $T$ is given by Bhandari and Khetan in \cite{Bhandari2024ImprovedUB}.

A natural restriction is to require the code $\C$ to be a linear code, that is a vector subspace of $\mathbb{F}_{3}^{n}$. With this restriction, one can examine the quantity $T_{L}(n)$, denoting the largest size of a linear trifferent code in $\mathbb{F}_{3}^{n}$. In \cite{Bishnoi_trifference} the authors show that a linear code is trifferent if and only if it is minimal, providing various significant improvements to existing knowledge on trifferent codes. In particular, they provide the best-known explicit construction of trifferent codes (which is only a constant factor away from the asymptotic lower bound on $T$). These works give excellent motivation for the study of minimal codes and strong blocking sets in the particular case when $q=3$.

\subsection{Open problems on minimal codes}

There are many ways to develop the theory of minimal codes.

\smallskip

First of all, since we have purposely defined minimal codes in as general a way as possible, it makes sense to examine minimal codes for various definitions of support. This has been the approach taken in the sum-rank metric by Santonostaso and Zullo in \cite{Subspacedesigns}, and by Borello and Zullo in \cite{sumrankminimal}. Since there are many more metrics than simply the Hamming, rank and sum-rank metrics, there is a lot of room for investigating different notions of support.

\smallskip

Since the theory of strong blocking sets is still relatively recent, as they were introduced in 2011, there still a lot to find out on them. In particular, stronger bounds on the minimal size of a strong blocking set in $\PG(k-1, q)$, noted by the function $m(k, q)$, would be a significant achievement.

Similarly, most upper bounds on $m(k, q)$ are probabilistic and do not provide explicit constructions. The best existing explicit constructions of small strong blocking sets, while asymptotically proportional to the bounds for $m(k, q)$, are quite large compared to $m(k, q)$. An interesting direction for further research is to try to find better explicit constructions, or constructions that can be implemented with small computational complexity. Since strong blocking sets have many interesting geometric properties, it seems that many different approaches should be possible. So far, many efficient constructions consist of an union of lines, but this still gives a lot of freedom for how to choose the lines in a suitable way. It is of course not required to consider sets of lines in order to obtain a strong blocking set.

\smallskip

In the rank metric, as explained in Section \ref{sec:rank}, there are only $k-1$ values of $n$ for which it is not known in general wether minimal codes exist. Settling wether minimal codes can exist for these values would help a lot in clarifying the general picture, as would providing efficient explicit constructions.

\smallskip

In Section \ref{sec:hamming} we have noted that the weight distribution of a minimal code must be quite restricted. This is also true in the rank metric, as noted in \cite{rankminimal}. In the opposite direction, the  Ashikhmin-Barg condition establishes that when the weight spectrum is sufficiently restricted, the code must be minimal. Therefore it would be interesting to further examine the relation between weight spectrum and minimality.

\smallskip

Finally, to the best of our knowledge, the dual codes of minimal codes have not been investigated in the nonbinary case. In analogy to the theory of the $2$-wise Davenport constant, there may be a simple combinatorial characterization of the duals of minimal codes. It remains to be seen what kind of mathematics will emerge from this kind of investigations.

\bibliographystyle{abbrv}
\bibliography{references.bib}

\end{document}